\documentclass[12pt]{article}

\usepackage{amsfonts}
\usepackage{amsmath}
\usepackage{amssymb}
\usepackage{amscd}
\usepackage{amsthm}
\usepackage{indentfirst}

\newtheorem{thm}{Theorem}[section]
\newtheorem{lemma}[thm]{Lemma}

\newtheorem{rem}[thm]{Remark}

\newtheorem{problem}[thm]{Problem}
\newtheorem{question}[thm]{Question}

\newcommand{\R}{{\mathbb{R}}}

\newcommand{\T}{{\mathbb{T}}}

\newcommand{\N}{{\mathbb{N}}}

\newcommand{\cF}{{\mathcal{F}}}

\def\id{{1\hskip-2.5pt{\rm l}}}

\newcommand{\tf}{{\widetilde{f}}}
\newcommand{\tg}{{\widetilde{g}}}

\newcommand{\tphi}{{\widetilde{\phi}}}
\newcommand{\tpsi}{{\widetilde{\psi}}}

\newcommand{\Hamc}{{{\it Ham^c}}}

\newcommand{\tHamc}{\widetilde{{\it Ham^c} }}

\newcommand{\Qed}{\hfill \qedsymbol \medskip}

\begin{document}

\title{$C^0$-rigidity of Poisson brackets \\
%{\it (preliminary version)}
}

\renewcommand{\thefootnote}{\alph{footnote}}

\author{\textsc Michael Entov$^{a}$,\ Leonid
Polterovich$^{b}$}

\footnotetext[1]{Partially supported by E. and J. Bishop Research
Fund and by the Israel Science Foundation grant $\#$ 881/06.}
\footnotetext[2]{Partially supported by the Israel Science
Foundation grant $\#$ 509/07.}

\date{\today}

\maketitle

\begin{abstract}
\noindent Consider a functional associating to a pair of compactly
supported smooth functions on a symplectic manifold the maximum of
their Poisson bracket. We show that this functional is lower
semi-continuous with respect to the product uniform ($C^0$) norm
on the space of pairs of such functions. This extends previous
results of Cardin-Viterbo and Zapolsky. The proof involves theory
of geodesics of the Hofer metric on the group of Hamiltonian
diffeomorphisms. We also discuss a failure of a similar
semi-continuity phenomenon for multiple Poisson brackets of three
or more functions.

\end{abstract}

%\tableofcontents

\renewcommand{\thefootnote}{\arabic{footnote}}
%\addtocounter{footnote}{-1}
%\vfil \eject

\section{Statement of results}
\label{sec-intro}

The subject of this note is function theory on symplectic
manifolds.  Let $(M, \omega)$ be a symplectic manifold (open or
closed). Denote by $C^\infty_c (M)$ the space of smooth compactly
supported functions on $M$ and by $\|\cdot \|$ the standard {\it
uniform norm}  (also called the {\it $C^0$-norm}) on it: $ \| F\|
:= \max_{x\in M} |F(x)|$.

The definition of the Poisson bracket $\{ F, G\}$ of two smooth
functions $F, G\in C^\infty_c (M)$ involves first derivatives of
the functions. Thus \emph{a priori} there is no restriction on
possible changes of $\{ F, G\}$ when $F$ and $G$ are slightly
perturbed in the uniform norm. Amazingly such restrictions do
exist: this was first pointed out  by F.Cardin and C.Viterbo
\cite{Car-Vit} who showed that
\[ \{ F, G\} \not\equiv 0
\Longrightarrow \liminf_{F',G'\stackrel{C^0}{\longrightarrow}F,G}
\|\{F',G'\}\|
>0.
\]
\noindent Our main result is as follows:

\begin{thm}
\label{thm-main} \begin{equation}\label{eq-thm-main}  \max \{ F, G\}
= \liminf_{F',G'\stackrel{C^0}{\longrightarrow}F,G} \max\{F',G'\}\
\end{equation}
for any symplectic manifold $M$ and any pair $F,G\in C^\infty_c
(M)$.

\end{thm}

\medskip
\noindent Replacing $F$ by $-F$, we get a similar result for
$-\min\{F,G\}$. In particular, this yields
\begin{equation}
\label{eq-zap}  \| \{ F, G\} \| =
\liminf_{F',G'\stackrel{C^0}{\longrightarrow}F,G} \|\{F',G'\}\|\;,
\end{equation}
which should be considered as a refinement of the Cardin-Viterbo
theorem, and which gives a positive answer to a question posed in
\cite{EPZ}.

In the case ${\rm dim}\, M=2$ formula \eqref{eq-zap} was first
proved by F.Zapolsky \cite{Z} by methods of two-dimensional
topology.

A generalization of the Cardin-Viterbo result in a different
direction has been found by V.Humili\`ere \cite{Hum}.

\begin{rem}
{\rm Note that \eqref{eq-zap} does {\it not} imply that $\{F',G'\}
\stackrel{C^0}{\longrightarrow} \{ F,G\}$ when
$F',G'\stackrel{C^0}{\longrightarrow} F,G $ -- see e.g. \cite{Hum}
for counterexamples. }
\end{rem}

\begin{rem}
{\rm Clearly $\liminf$ cannot be replaced in the theorem by $\lim$:
the maximum of the Poisson bracket of two functions can be
arbitrarily increased by arbitrarily $C^0$-small perturbations of
the functions.}
\end{rem}

In the proof of Theorem~\ref{thm-main} we use the following
ingredient from ``hard" symplectic topology: Denote by $\Hamc (M)$
the group of Hamiltonian diffeomorphisms of $M$ generated by
Hamiltonian flows with compact support. Then sufficiently small
segments of one-parameter subgroups of the group $\Hamc (M)$ of
Hamiltonian diffeomorphisms of $M$ minimize the ``positive part of
the Hofer length" among all paths on the group in their homotopy
class with fixed end points. This was proved  by D.McDuff in
\cite[Proposition 1.5]{McD-variants} for closed manifolds and in
\cite[Proposition 1.7]{McD-monodromy} for open ones; see also
papers \cite{Bialy-Pol}, \cite{Lal-McD}, \cite{En},
\cite{McD-Slim}, \cite{KL} \cite{Oh} for related results in this
direction.

After an early draft of this paper has been written, L.Buhovsky
found a different proof of Theorem \ref{thm-main} based on an
ingenious application of the energy-capacity inequality. Buhovsky's
method enables him to give a quantitative estimate on the rate of
convergence in the right-hand side of \eqref{eq-thm-main}. These
results will appear in a forthcoming article \cite{Buh}.

The next result gives an evidence for a failure of $C^0$-rigidity
for multiple Poisson brackets.

\begin{thm}
\label{thm-multiple-brackets} Let $M$ be a symplectic manifold.
There exists a constant $N\in \N$, depending only on the dimension
of $M$, such that {\bf for any} smooth functions $F_1,\ldots,
F_N\in C^\infty_c (M)$ there exist $F'_1,\ldots, F'_N \in
C^\infty_c (M)$ arbitrarily close in the uniform norm,
respectively, to $F_1,\ldots, F_N$  which satisfy the following
relation:
\[
\{ F'_1, \{ F'_2, \ldots \{F'_{N-1}, F'_N\} \}\ldots \} \equiv 0.
\]
\end{thm}

\medskip
\noindent We shall see in Section~\ref{sec-multbr} below that in
the case $\dim M=2$ the result above holds for $N=3$.

\begin{question}\label{q-1}{\rm Does the theorem above remain valid
with $N=3$ on an arbitrary symplectic manifold?}
\end{question}

\medskip
\noindent The following claim, though it does not answer
Question~\ref{q-1}, shows that Theorem~\ref{thm-main} cannot be
formally extended to the triple Poisson bracket.

\begin{thm}
\label{thm-3-brackets} For any symplectic manifold $M$ {\bf one
can find} 3 functions $F,G,H\in C^\infty_c (M)$ satisfying $\{F,\{
G,H\}\}\not\equiv 0$ such that there exist smooth functions $F',
G', H' \in C^\infty_c (M)$ arbitrarily close in the uniform norm,
respectively, to $F, G,H$ and satisfying the condition
$$\{ F',\{G',H'\}\}\equiv 0.$$
\end{thm}

\medskip
\noindent The theorem will be proved in Section~\ref{sec-multbr}.
The proof shows that the phenomenon is local: we just implant a
2-dimen\-si\-onal example (see the remark after
Theorem~\ref{thm-multiple-brackets}) in a Darboux chart.

Surprisingly, the next problem is open even in dimension 2:

\begin{problem}\label{prob}{\rm
Compare
$$\liminf_{{F',G'\stackrel{C^0}{\longrightarrow}F,G}}\max\{\{F',G'\},G'\}$$
with $\max \{\{F,G\},G\}$ for some/all pairs of functions $F,G$ on
some/all symplectic manifolds.}
\end{problem}

\section{Proofs}
\label{sec-proofs}

\subsection{Preliminaries}

Given a (time-dependent) Hamiltonian $H: M\times [0,1]\to \R$,
denote by $X_H$ the (time-dependent) Hamiltonian vector field
generated by $H$. The Poisson bracket of two functions $F,G \in
C^{\infty}_c (M)$ is defined by $\{F,G\}= dF(X_G)$.

Let $\Hamc (M)$ be the group of Hamiltonian diffeomorphisms of
$(M,\omega)$ generated by compactly supported (time-dependent)
Hamiltonians. Write $\tHamc (M)$ for the universal cover of $\Hamc
(M)$, where the base point is chosen to be the identity map $\id$.
Denote by $\psi_H^t$, $t\in\R$, the Hamiltonian flow generated by
$H$ (i.e. the flow of $X_H$). Let $\psi_H := \psi_H^1$ and let
$\tpsi_H\in \tHamc (M)$ be the lift of $\psi_H$ associated to the
path $\{\psi_H^t\}, t \in [0;1]$. We will say that $\psi_H$ and
$\tpsi_H$ {\it are generated by} $H$. We will also denote $\| H\|
= \max_{M\times [0,1]} | H(x,t)|$ (for time-independent
Hamiltonians this norm coincides with the uniform norm on
$C^\infty_c (M)$ introduced above). Set $H_t = H(\cdot, t)$.

Recall that the flow $\psi_H^t\psi_K^t$ is generated by the
Hamiltonian $H\sharp K (x,t) = H (x,t) + K((\psi_H^t)^{-1} x, t)$
and the flow $\psi_H \psi_K^t (\psi_H)^{-1}$  by $K \big(
(\psi_H)^{-1} x, t\big)$.

A Hamiltonian $H$ on $M \times [0;1]$ is called {\it normalized} if
either $M$ is open and $\bigcup_t \text{support}(H_t)$ is contained
in a compact subset of $M$, or $M$ is closed and $H_t$ has zero mean
for all $t$. The set of all normalized Hamiltonian functions is
denoted by $\cF$. Note that if $H,K\in\cF$ then both $H\sharp K$ and
$K \big( (\psi_H)^{-1} x, t \big)$ also belong to $\cF$.

For $a,b \in \tHamc (M,\omega)$ write $[a,b]$ for the commutator
$aba^{-1}b^{-1}$.

\begin{lemma}
\label{lem-1} Assume $H,K\in C^\infty_c (M)$ are time-independent
Hamiltonians. Then $[\tpsi_H, \tpsi_K]$ can be generated by $L
(x,t) = H (x) - H \big( \psi_K^{-1}\psi_H^{-t} x\big)$.

\end{lemma}

\begin{proof} It is easy to see that the element $[\tpsi_H, \tpsi_K] \in
\tHamc (M)$ can be represented by the path $\{\psi_H^t \psi_K
\psi_H^{-t}\psi_K^{-1}\}$ where $t\in [0;1]$. The flow
$\psi_H^{-t}$ is generated by $-H$ and therefore the flow $\psi_K
\psi_H^{-t} \psi_K^{-1}$ is generated by the Hamiltonian $-H\circ
\psi_K^{-1}$. Thus the flow $\psi_H^t \psi_K \psi_H^{-t}
\psi_K^{-1}$ is generated by $H\sharp \big( -H\circ \psi_K^{-1}
\big) (x,t) = H (x) - H \big( \psi_K^{-1}\psi_H^{-t} x\big)$.
\end{proof}

The group $\tHamc (M)$ carries conjugation-invariant functionals
$\rho^+$ and $\rho$ defined by
\[\rho^+(\tpsi) :=
\inf_H \int_0^1 \max_{x \in M} H(x,t)\;dt
\]
and
\[\rho(\tpsi) :=
\inf_H \int_0^1 (\max_{x \in M} H(x,t) - \min_{x \in M}
H(x,t))\;dt\;,
\]
where the infimum is taken over all Hamiltonians $H\in\cF$
generating $\tpsi$. The functional $\rho$ is the Hofer (semi)-norm
\cite{H} (see e.g. \cite{Pol-book} for an introduction to Hofer's
geometry). It gives rise to the bi-invariant Hofer (pseudo-)metric
on $\tHamc (M)$ by  $d(\tphi,\tpsi) = \rho(\tphi^{-1}\tpsi)$. The
functional $\rho^+$, which is sometimes called the ``positive part
of the Hofer norm", satisfies the triangle inequality but is not
symmetric. Note also that $\rho^+ \leq \rho$. We shall use the
following properties of these functionals. By the triangle
inequality for $\rho^+$
\begin{equation}\label{eqn-hof-1}
|\rho^+(\tphi)-\rho^+(\tpsi)| \leq
\max(\rho^+(\tphi^{-1}\tpsi),\rho^+(\tpsi^{-1}\tphi)) \leq
d(\tphi,\tpsi)\;.
\end{equation}
This readily yields
\begin{equation}\label{eqn-hof-2}
|\rho^+(\tpsi_H)-\rho^+(\tpsi_K)| \leq d(\tpsi_H,\tpsi_K) \leq
2||H-K||
\end{equation}
for any $H,K\in\cF$. McDuff showed \cite[Proposition
1.5]{McD-variants}, \cite[Proposition 1.7]{McD-monodromy} that for
every {\it time-independent} function $H \in \cF$ there exists
$\delta > 0$ so that
\begin{equation}\label{eq-geodes}
\rho^+(\tpsi_{tH}) = t\cdot \max H  \;\;\forall t \in
(0;\delta)\;.
\end{equation}

\begin{lemma}
\label{lem-2} Assume $H,K\in C^\infty_c (M)$ are time-independent
Hamiltonians with zero mean. Then $\rho^+([\tpsi_H,\tpsi_K]) \leq
\max\{ H,K\}$.
\end{lemma}

\begin{proof}
By Lemma~\ref{lem-1} $[\tpsi_H,\tpsi_K]$ can be generated by
\[L (x,t) = H (x) - H \big( \psi_K^{-1}\psi_H^{-t} x\big).
\] Note that
\[
\int_0^1\max L(x,t)\;dt  = \int_0^1 \max ( H  - H \circ
\psi_K^{-1} \circ \psi_H^{-t})\; dt\]\[ = \int_0^1 \max( H \circ
\psi_H^t - H \circ \psi_K^{-1})\;dt
\] \[ =  \int_0^1\max ( H - H \circ \psi_K^{-1})\;dt = \int_0^1 \max (H \circ \psi_K  -
H)\;dt
\]
since $H$ is constant on the orbits of the flow $\psi_H^t$. Taking
into account that
\[ H (\psi_K x) - H (x) = \int_0^1 \frac{d}{dt} H
(\psi_K^t x) dt  = \int_0^1 \{ H, K\} (\psi_K^t x) dt,\] we get
that
\[
\rho^+([\tpsi_H,\tpsi_K])\leq \int_0^1\max L (x,t)\;dt  \leq \max
\{H, K\},
\]
which yields the lemma.
\end{proof}

\subsection{Proof of Theorem~\ref{thm-main}}

We assume without loss of generality that all the functions $F_i,
G_i, F,G$ are normalized.

Denote by $f_s$, $g_t$, $s,t\in[0,1]$, the Hamiltonian flows
generated by $F$ and $G$, and by $\tf_s,\tf_t$ their respective
lifts to $\tHamc (M)$. Note that for fixed $s$ and $t$ the
elements $\tf_s$ and $\tg_t$ are generated, respectively, by the
Hamiltonians $sF$ and $tG$.

By Lemma~\ref{lem-1} for fixed $s,t$ the commutator $[\tf_s,
\tg_t]$ can be generated by the Hamiltonian $L_{s,t} (x,\tau) = sF
(x) - sF (g_t^{-1} f_{\tau s}^{-1} x)$ (use Lem\-ma~\ref{lem-1}
with $H=sF$, $K=tG$ and note that $\psi_{sF}^\tau = f_{\tau s}$).
Clearly $L_{s,t}\in\cF$ since $F,G\in \cF$.

\begin{lemma}
\label{lem-3}

$L_{s,t} = st\{ F,G\} + K_{s,t}$, where $\| K_{s,t} \|/st \to 0$
as $s,t\to 0$.

\end{lemma}

\begin{proof}
We need to compute the relevant terms in the expansion of
$L_{s,t}$ with respect to $s,t$ at $s=0, t=0$.

Clearly, $L_{0,0} \equiv 0$.

The first order terms are as follows:
\[
\frac{\partial L_{s,t}}{\partial s}  (x,\tau) =
\partial \big(sF (x) - sF (g_t^{-1} f_{\tau s}^{-1} x) \big) /\partial
s =
\]
\[
= F(x) - F (g_t^{-1} f_{\tau s}^{-1} x) - s dF\circ dg_t^{-1}
(X_{-sF} (x) ) =
\]
\[
= F(x) - F (g_t^{-1} f_{\tau s}^{-1} x) + s^2 dF\circ dg_t^{-1}
(X_F (x) ),
\]
and
\[
\frac{\partial L_{s,t}}{\partial t}  (x,\tau) =
\partial \big(sF (x) - sF (g_t^{-1} f_{\tau s}^{-1} x) \big) /\partial
t = \] \[ = - sdF \big( X_{-G} (f_{\tau s}^{-1} x) \big) = s\{
F,G\} (f_{\tau s}^{-1} x).
\]
Evaluating ${\partial L_{s,t}}/{\partial s} (x,\tau)$ and
${\partial L_{s,t}}/{\partial t} (x,\tau)$, respectively, at the
points $(s,0)$ and $(0,t)$ (for a fixed $(x,\tau)$) we see that
\[
\frac{\partial L_{s,0}}{\partial s} (x,\tau) \equiv 0 \] (since
$F$ is constant on the orbits of the flow $f_s$, $s\in \R$) and
\[
\frac{\partial L_{0,t}}{\partial t} (x,\tau) \equiv 0. \] Thus
\[
\left. \frac{\partial^k}{\partial s^k} \right|_{(s,t)=(0,0)}
L_{s,t} (x,\tau) = 0  = \left. \frac{\partial^k}{\partial t^k}
\right|_{(s,t)=(0,0)} L_{s,t} (x,\tau), \ {\rm for\ any}\ k\geq 1.
\]

Finally, let us compute $\left. \frac{\partial^2 }{\partial s
\partial t} \right|_{(s,t)=(0,0)} L_{s,t} (x,\tau)$:
\[
\left. \frac{\partial^2}{\partial s \partial t}
\right|_{(s,t)=(0,0)} L_{s,t} (x,\tau) = \left.
\frac{\partial}{\partial s} \right|_{s=0} \frac{\partial
L_{s,0}}{\partial t} (x,\tau) = \]
\[
= \left. \frac{\partial}{\partial s} \right|_{s=0} s\{ F,G\}
(f_{\tau s}^{-1} x) = \{ F,G\} (x).
\]
This finishes the proof of the lemma.
\end{proof}

Now we are ready to complete the proof of Theorem~\ref{thm-main}.
The inequality \[\max \{ F, G\} \geq
\liminf_{F',G'\stackrel{C^0}{\longrightarrow}F,G} \max \{F',G'\}
\] is trivial so we only need to prove the opposite one. Let $F_i$,
$G_i$ be sequences of smooth functions such that
\[F_i, G_i \stackrel{C^0}{\longrightarrow} F,G, \ i\to +\infty, \] and
\[\max \{ F_i,
G_i\} \to A, \ i\to +\infty.\]  We need to show that $\max \{ F, G\}
\leq A$.

Assume on the contrary that $\max \{ F, G\}  > A$. Pick $B$ such
that $A<B< \max \{ F, G\} $. Then for any sufficiently large $i$
\[
\max \{ F_i, G_i\} \leq B.
\]
Denote by $f_{s,i}, g_{t,i}$, respectively, the time-$s$ and
time-$t$ maps of the flows generated by $F_i$ and $G_i$. Their
lifts to $\tHamc (M)$ will be decorated by tildes. The right
inequality in (\ref{eqn-hof-2}) easily implies  that the sequences
$\tf_{s,i}$ and $\tg_{t,i}$ converge, respectively, to $\tf_s$ and
$\tg_s$ in the Hofer (pseudo-)metric. Since by \eqref{eqn-hof-1}
the functional $\rho^+$ is continuous in the Hofer
(pseudo-)metric,
\[
\rho^+ ([\tf_{i,s}, \tg_{i,t}])\to  \rho^+ ([\tf_{s}, \tg_{t}]) \;\;
\text{as}\;\; i \to \infty\;.
\]
By Lemma~\ref{lem-2},
\[
\rho^+ ([\tf_{i,s}, \tg_{i,t}])\leq st \cdot \max \{ F_i, G_i\} \leq
stB
\]
for any sufficiently large $i$. Hence, taking the limit in the
left-hand side as $i\to +\infty$, we get
\begin{equation}
\label{eqn-Bst} \rho^+ ([\tf_s, \tg_t])\leq stB.
\end{equation}
Choose $\epsilon > 0$ such that $B+2\epsilon< \max \{ F, G\} $. Take
sufficiently small $s,t>0$ so that the function $K_{s,t}$ from
Lemma~\ref{lem-3} admits a bound
\begin{equation}
\label{eqn-K-s-t} \| K_{s,t}\| \leq \epsilon st
\end{equation}
 and so that the Hamiltonian $st\{F,G\}$ is sufficiently small and
satisfies
\begin{equation}
\label{eqn-C0-flatness} \rho^+ (\tpsi_{st\{F,G\}}) = st\cdot\max
\{ F, G\}\;,
\end{equation}
see formula \eqref{eq-geodes}. Lemma~\ref{lem-3} and  inequalities
(\ref{eqn-K-s-t}), (\ref{eqn-hof-2}) yield
\[
|\rho^+ ([\tf_s, \tg_t])-\rho^+(\tpsi_{st\{F,G\}})| \leq 2\epsilon
st\;.
\]
Hence,
\[
\rho^+ ([\tf_s, \tg_t]\geq \rho^+ (\tpsi_{st\{F,G\}}) - 2\epsilon
st = st (\max \{ F, G\}  - 2\epsilon).
\]
Combining this with (\ref{eqn-Bst}), we get
\[
st (\max \{ F, G\}  - 2\epsilon) \leq \rho^+ ([\tf_s, \tg_t])\leq
stB,
\]
and hence
\[
 \max \{ F, G\}  - 2\epsilon \leq B
\]
which contradicts our choice of $B$ and $\epsilon$. We have obtained
a contradiction. Hence $\max \{ F, G\}  \leq A$ and the theorem is
proven.\Qed

\subsection{Proofs of
Theorems~\ref{thm-multiple-brackets},
\ref{thm-3-brackets}}\label{sec-multbr}

\medskip \noindent
{\bf Proof of Theorem~\ref{thm-multiple-brackets}.}

For simplicity we will prove the result in the case ${\rm dim}\, M
= \T^2$ with $N=3$. The general case can be done in a similar way
using \cite{Rudyak-Schlenk}.

Define a {\it thick grid} $T$ with mesh $c$ in $M$ as a union of
pair-wise disjoint squares on $M$ such that each square has a side
$2c$ and the centers of the squares form a rectangular grid with
the mesh $3c$. A {\it $T$-tamed function} is a smooth function
which is constant in a small neighborhood of each square of the
thick grid $T$ (but its values may vary from square to square).

One can easily construct a sequence $c_i \to 0$ and $N=3$ thick
grids $U_i,V_i,W_i$ with mesh $c_i$ so that $U_i\cup V_i \cup W_i
= M$ for all $i$. (See \cite{Rudyak-Schlenk} on how to construct a
similar covering of an arbitrary $M$ by a number of thick grids
depending only on ${\rm dim}\, M$).

Now for every $\epsilon > 0$ there exists $i$ large enough so that
every triple of functions $F_1,F_2,F_3\in C^\infty_c (M)$ can be
$\epsilon$-approximated, respectively, by $U_i,V_i,W_i$-tamed
functions $F'_1,F'_2,F'_3 \in C^\infty_c (M)$. Take any point $x
\in M$. Then at least one of the functions $F'_1,F'_2,F'_3$ is
constant near $x$. Thus $\{F'_1,\{F'_2,F'_3\}\} \equiv 0$, and the
claim follows. \Qed

\bigskip
\noindent
\medskip \noindent
{\bf Proof of Theorem~\ref{thm-3-brackets}.}

Assume ${\rm dim}\, M = 2n > 2$ (the case ${\rm dim}\, M =2$ has
been dealt with in the proof of
Theorem~\ref{thm-multiple-brackets}). In a local Darboux chart
with coordinates $p_1, q_1,\ldots, p_n, q_n$ on $M$ choose an open
cube $$P= K^{2n-2} \times K^2,$$ where $K^{2n-2}$ is an open cube
in the $(p_1, q_1,\ldots, p_{n-1}, q_{n-1})$-coordinate plane and
$K^2$ is a open square in the $(p_n, q_n)$-coordinate plane. Fix a
smooth compactly supported non-zero function $\chi$ on $K^{2n-2}$.
Given a smooth compactly supported function $L$ on $K^2$, define
the function $\chi L\in C^\infty_c (M)$ as
$$\chi L (p_1, q_1,\ldots, p_n, q_n) := \chi (p_1, q_1,\ldots, p_{n-1}, q_{n-1}) L (p_n, q_n)$$
on $P$ and as zero outside $P$.

Now pick any functions $F_1,G_1,H_1\in C^\infty_c (K^2)$ such that
$$\{ F_1,\{ G_1,H_1\}\}\not\equiv 0.$$
Set $$ F:= \chi F_1, G:=\chi G_1, H:= \chi H_1 \in C^\infty_c
(M).$$ As in the proof of Theorem~\ref{thm-multiple-brackets}
(note that in the case of the two-dimensio\-nal square the
construction of the thick grids is as easy as in the case of
$\T^2$), choose $C^0$-small perturbations $F'_1, G'_1, H'_1\in
C^\infty_c (K^2)$ of $F_1,G_1,H_1$ so that
$$\{ F'_1,\{ G'_1,H'_1\}\} \equiv 0.$$
Then $F':=\chi F'_1, G':=\chi G'_1, H':=\chi H'_1\in C^\infty_c
(M)$ satisfy
$$\{ F', \{ G', H'\}\} = \{ \chi F'_1,\{ \chi G'_1, \chi H'_1\}\}
=
\chi^3 \{ F'_1,\{ G'_1,  H'_1\}\} \equiv 0,$$ because of the
Leibniz rule for Poisson brackets and because the Poisson bracket
of $\chi$ and any function of $p_n, q_n$ vanishes identically. For
the same reason
$$\{ F,\{ G, H\}\} = \{ \chi F_1, \{ \chi G_1, \chi H_1\}\} = \chi^3
\{ F_1,\{ G_1,H_1\}\} \not\equiv 0.$$ Clearly, by choosing $F'_1,
G'_1, H'_1$ arbitrarily $C^0$-close to $F_1,G_1,H_1$ in
$C^\infty_c (K^2)$ we can turn $F', G', H'$ into arbitrarily
$C^0$-small perturbations of $F, G, H$ in $C^\infty_c (M)$. Thus
we have constructed $F,G,H, F',G',H'$ satisfying the required
conditions. \Qed

\medskip
\noindent {\bf Acknowledgement.} We thank L.Buhovsky, D.Burago,
M.Khanevsky and D.McDuff for useful discussions and comments.

\medskip

\bibliographystyle{alpha}

\bigskip

\setlength{\parindent}{-0.14in}

\begin{tabular}{ll}
Michael Entov & Leonid Polterovich\\
Department of Mathematics&School of Mathematical Sciences\\
Technion - Israel Institute of Technology & Tel Aviv University\\
Haifa 32000, Israel & Tel Aviv 69978, Israel\\
entov@math.technion.ac.il & polterov@post.tau.ac.il\\
\end{tabular}

\end{document}